\documentstyle[12pt]{article}
\newcommand{\too}{\longrightarrow}
\newcommand{\om}{\omega}
\newcommand{\an}{{\pi_{\#}}}

\newcommand{\U}{{\cal U}}
\newcommand{\X}{{\cal X}}
\newcommand{\G}{{\cal G}}
\newcommand{\F}{{\cal F}}

\newcommand{\M}{{\cal M}}

\newcommand{\D}{{\cal D}}

\newcommand{\di}{\displaystyle}
\newcommand{\Om}{\Omega}
\newcommand{\na}{\nabla}
\newcommand{\wi}{\widetilde}
\newcommand{\al}{\alpha}
\newcommand{\be}{\beta}
\newcommand{\ga}{\gamma}
\newcommand{\Ga}{\Gamma}
\newcommand{\e}{\epsilon}
\newcommand{\la}{\lambda}

\newcommand{\inj}{\hookrightarrow}

\def \reel{ {\rm I}\!{\rm R} }

 \def \comp{ \;{}^{ {}_\vert }\!\!\!{\rm C}   }

 \def \rat{ {\rm Q}\kern-.65em {}^{{}_/ }}

\newtheorem{th}{Theorem}[section]
\newtheorem{pr}{Proposition}[section]

\title{Solutions of the classical Yang-Baxter equation and noncommutative
deformations
 } \author{M. Boucetta\footnote{Recherche men\'ee dans le cadre du Programme
 Th\'ematique d'Appui \`a la Recherche Scientifique PROTARS III.}} \date{ }
\begin{document}
\maketitle
\parindent=0cm

{\bf Abstract.}
 In this paper, I will show that, if a Lie algebra $\G$  acts on a
 manifold $P$, any solution of the classical
 Yang-Baxter equation on  $\G$ gives arise to a Poisson tensor  on $P$ and a torsion-free
 and flat
 contravariant connection (with respect to the Poisson tensor).
 Moreover, if the action is locally free,  the matacurvature of
 the above contravariant connection vanishes. This
 will permit to get a large class of manifolds which satisfy the
 necessary conditions, presented by Hawkins in
 arXiv:math.QA/050423, to the existence of a noncommutative
 deformation.
\bigskip

\section{ Introduction and main results}
In [9], Hawkins showed that  a deformation of the differential
graded algebra of differential forms $\Om^*(P)$ on a manifold $P$
gives arise to a Poisson tensor on $P$ (which characterizes the
deformation) and a torsion-free and flat contravariant connection
whose metacurvature vanishes. In Riemannian case, Hawkins showed
that if a Riemannian manifold $P$ is deformed into a real spectral
triple, the Riemannian metric and the Poisson tensor (which
characterizes the deformation) are compatible in the following
sense:
\begin{enumerate}
\item the metric contravariant connection $\D$ associated to the
metric and the Poisson tensor is flat. \item The metacurvature of
$\D$ vanishes. \item The Poisson tensor $\pi$ is compatible with
the Riemannian volume $\e$:
$$d(i_{\pi}\e)=0.$$\end{enumerate}

The metric contravariant connection is a torsion-free
contravariant connection associated naturally to any couple of
Riemannian metric and Poisson tensor. It has appeared first in
[3]. The metacurvature, introduced by Hawkins in [9], is a
$(2,3)$-tensor field (symmetric in the contravariant indices and
antisymmetric in the covariant indices) associated naturally to
any torsion-free and flat contravariant connection.

In this paper, we will construct a large class of smooth manifolds
which satisfy the necessary conditions (presented by Hawkins in
[9]) to deform the differential graded algebra of differential
forms. We will also give a large class of (pseudo)-Riemannian
manifolds with a Poisson tensor compatible with the metric in the
sense of Hawkins.

In order to state the main results of this paper, let us recall
some classical results on symplectic Lie groups and on  solutions
of the classical Yang-Baxter equation.

Let $G$ be a Lie group and $\G$ its Lie algebra. The group $G$ is
a sympelctic Lie group if there exists on $G$ a left invariant
symplectic form. It is clear that $G$ is a symplectic Lie group if
and only if there exists on $\G$ a non-degenerate 2-form $\om$
such that
$$\om([x,y],z)+\om([y,z],x)+\om([z,x],y)=0,\quad
x,y,z\in\G.\eqno(1).$$

A Lie algebra with a 2-form such above is called a symplectic Lie
algebra. It is well known since
 the end of 1960s (see e.g. [5]) that, if $(\G,\om)$ is a symplectic Lie algebra,
  the formula
 $$\om(A_xy,z)=-\om(y,[x,z])\eqno(2)$$defines a product
 $A:\G\times \G\too \G$ which verifies:\begin{enumerate}
 \item $A_xy-A_yx=[x,y];$\item
 $A_{[x,y]}z-(A_xA_yz-A_yA_xz)=0.$\end{enumerate}Thus $A$ defines
 on a $G$  a flat  and
 torsion-free linear (covariant) connection (or equivalently an affine
 structure) which is left invariant.

Let $\G$ be a Lie algebra and  $r\in\G\wedge\G$. We will also
denote by $r:\G^*\too\G$ the linear map induced by $r$. The
bi-vector $r$ satisfies the classical Yang-Baxter equation if
$$[r,r]=0,\eqno(Y-B)$$where $[r,r]\in\G\wedge\G\wedge\G$ is defined
by
$$[r,r](\al,\be,\ga)=\al([r(\be),r(\ga)])+\be([r(\ga),r(\al)])+
\ga([r(\al),r(\be)]).$$

Solutions of $(Y-B)$ are strongly related to symplectic Lie
algebras and hence to left invariant affine structures. Let us
explain this relation.

One can observe that to give $r\in {\cal G}\wedge{\cal G}$ is
equivalent to give a vectorial subspace $S_r\subset\cal G$ and a
non-degenerate 2-form $\om_r\in\wedge^2S^*_r$. Indeed, for $r\in
{\cal G}\wedge{\cal G}$, we put $S_r=Imr$ and
$\om_r(u,v)=r(r^{-1}(u),r^{-1}(v))$ where $u,v\in S_r$ and
$r^{-1}(u)$ is any antecedent of $u$ by $r$.

Conversely, let $(S,\om)$ be  a vectorial subspace of $\cal G$
with a non-degenerate 2-form. The 2-form $\om$ defines an
isomorphism $\om^b:S\too S^*$ by $\om^b(u)=\om(u,.)$, we denote by
$\om^\#:S^*\too S$ its inverse and we put $r=\om^\#\circ i^*$
where $i^*:{\cal G}^*\too S^*$ is the dual of the inclusion
$i:S\inj\cal G$.

With this observation in mind, the following classical proposition
gives another description of the solutions of $(Y-B)$.
\begin{pr} Let $r\in {\cal G}\wedge{\cal G}$ and $(Imr,\om_r)$ its
associated subspace. The following assertions are equivalent:
\begin{enumerate}\item
 $r$ is a solution of $(Y-B)$.
\item $(Imr,\om_r) $ is a symplectic subalgebra of $\cal
G$.\end{enumerate}
\end{pr}

A solution $r$ of $(Y-B)$ will be called abelian, unimodular
etc...  if $Imr$ is an abelian subalgebra, an unimodular
subalgebra etc...
\bigskip

 {\bf Remark.}\begin{enumerate}\item
  Let $\G$ be a  Lie algebra and $S$ an even dimensional abelian
subalgebra of $\G$. Any non-degenerate 2-form $\om$ on $S$
verifies the assertion  2. in Proposition 1.1 and hence $(S,\om)$
defines a solution of $(Y-B)$. \item Remark that any solution of
$(Y-B)$ on a compact Lie algebra is abelian.\end{enumerate}

 Let $P$ a smooth manifold and $\G$ a Lie algebra which acts on
 $P$, i.e., there exists a morphism of Lie algebras
 $\Ga:\G\too\X(P)$ from $\G$ to the Lie algebra of vector fields
 on $P$.
  Let $r\in\wedge^2\G$ be a solution of $(Y-B)$. If
$r=\sum_{i<j}a_{ij}u_i\wedge u_j$, put
$$\D^{\Ga}_\al\be:=\sum_{i<j}a_{ij}\left(\al(U_i)L_{U_j}\be-\al(U_j)L_{U_i}\be
\right),\eqno(3)$$where $\al,\be$ are two differential 1-forms on
$P$ and $U_i=\Ga(u_i)$.

  The restriction of $\Ga$ to $Imr$ defines on $P$ a
 singular foliation which coincides with the symplectic foliation
 of the Poisson tensor $\Ga(r)$ if the action of $Imr$ is locally free.
 In this case, since $Imr$ is a symplectic Lie algebra, the product
 given by $(2)$ defines, on any symplectic
 leaf of $\Ga(r)$, a flat and torsion-free covariant connection .
 Thus any symplectic leaf of $\Ga(r)$ becomes  an affine
 manifold. In fact, our main result in this paper, asserts that
 $\D^\Ga$ given by $(3)$ is a
 flat and torsion-free
 contravariant connection associated to $\Ga(r)$. Hence $\D^\Ga$  induces on
 any symplectic leaf an affine structure which is the affine structure described above.

 Let us
 state now our main results.
 \begin{th} Let $P$ be a differentiable manifold, $\G$ a Lie algebra
with  $\Ga:\G\too \X(P)$ a Lie algebras morphism from $\G$ to the
Lie algebra of  vector fields and $r\in\wedge^2\G$ a solution of
the classical Yang-Baxter equation such that $Imr$ is an
unimodular Lie algebra. Then, for any volume form $\e$ on $P$ such
that $L_{\Ga(u)}\e=0$ for each $u\in Imr$, $\Ga(r)$ is compatible
with $\e$, i.e., $$d(i_{\Ga(r)}\e)=0.$$\end{th}
\begin{th} Let $P$ be a differentiable manifold, $\G$ a Lie algebra
with  $\Ga:\G\too \X(P)$ a Lie algebras morphism from $\G$ to the
Lie algebra of  vector fields and $r\in\wedge^2\G$ a solution of
the classical Yang-Baxter equation. Then:

$i)$ $\D^\Ga$ given by (3) defines a torsion-free and flat
contravariant connection $\D^{\Ga}$ associated to the Poisson
tensor $\Ga(r)$ which depends only on $r$ and $\Ga$.

$ii)$ If $P$ is a Riemannian manifold and $\Ga(u)$ is a Killing
vector field for any $u\in Imr$, then $\D^{\Ga}$ is the metric
contravariant connection associated to the metric and $\Ga(r)$.

$iii)$ If the restriction of $\Ga$ to $Imr$ is a locally free
action then the metacurvature of $\D^\Ga$ vanishes.\end{th}

 There are some interesting implications of Theorem 1.1 and
 Theorem 1.2:
\begin{enumerate}
\item Let $G$ be a Lie group and $\G$ its  Lie algebra identified
to $T_eG$. Let $\Ga^l:\G\too\X(G)$ and $\Ga^r:\G\too\X(G)$ be
respectively the left and the right action of $\G$ on $G$.  For
any $\al\in\G^*$, we denote by $\al^l$ (resp. $\al^r$) the
left-invariant (resp. right-invariant) differential 1-form on $G$
associated to $\al$. For any $r\in\wedge^2\G$  a solution of
$(Y-B)$, one can check easily that $\D^{\Ga^l}$ and $\D^{\Ga^r}$
are given by
$$\D^{\Ga^l}_{\al^l}\be^l=-(ad^*_{r(\al)}\be)^l\quad
\mbox{and}\quad \D^{\Ga^r}_{\al^r}\be^r=(ad^*_{r(\al)}\be)^r.$$
Let  $<,>^l$ (resp. $<,>^r$) be a left-invariant (resp.
right-invariant) pseudo-Riemannian metric on $G$. If $r$ is
unimodular, then the Poisson tensor $\Ga^r(r)$ (resp. $\Ga^l(r)$)
is compatible in the sense of Hawkins with $<,>^l$ (resp.
$<,>^r$). \item Since any Lie algebra admits a non trivial
solution of $(Y-B)$ (see [6]), according to Theorem 1.2, any
locally free action of $\G$ on a manifold $P$ gives arise to a non
trivial Poisson tensor with a non trivial torsion-free and flat
contravariant connection whose metacurvature vanishes; so the
necessary conditions to deform the differential graded algebra of
differential forms $\Om^*(P)$ are satisfied (see [9]). \item Any
locally free action by isometries of an unimodular symplectic Lie
group $G$ on a pseudo-Riemannian manifold $P$ gives arise on $P$
to a Poisson tensor which is compatible with the metric in the
sense of Hawkins. \item Any symplectic nilmanifold is a quotient
of nilpotent  symplectic Lie group by a discrete co-compact
subgroup (see [1]), so any symplectic nilmanifold admit a
torsion-free and flat contravariant connection whose metacurvature
vanishes. Symplectic  nilpotent groups were classified by Medina
and Revoy [14]. \item The affine group $K^n\times GL(K^n)$
($K=\reel$ or $\comp$) admits many left invariant symplectic forms
(see [2]) which implies that its Lie algebra carries many
invertible solutions of $(Y-B)$. So any locally free action of the
affine group on a manifold gives arise to a Poisson tensor with a
non trivial, torsion-free and flat contravariant connection whose
metacurvature vanishes.\item Let $G$ be a Lie group with a
bi-invariant pseudo-Riemannian metric and $r$ an unimodular
solution of $(Y-B)$. For any discrete, co-compact
 subgroup $\Lambda$ of $G$,
 $G$ acts on the compact manifold
 $P:=G/\Lambda$ by isometries, so we get a Poisson tensor on the
 compact pseudo-Riemannian
 manifold
 $P$ compatible with the pseudo-Riemannian metric in the sense of
 Hawkins.

A connected Lie group $G$ admits a bi-invariant Riemannian metric
if and only if it is isomorphic to the cartesian product of a
compact group and a commutative group (see [14]). Any solution of
$(Y-B)$ on the Lie algebra of a such group is abelian. So any
triple $(G,\Lambda,S)$, where $G$ is a Lie group with a
bi-invariant Riemannian metric, $\Lambda$ is a discrete and
co-compact
 subgroup  of $G$ and $S$ an even dimensional subalgebra of the
 Lie algebra of $G$, gives arise to a compact Riemannian manifold
 with a Poisson tensor compatible with the metric in the sense of
 Hawkins.

Connected Lie groups which admit a bi-invariant pseudo-Riemannian
metric were classified in [11]. We will give now an example of a
compact Lorentzian manifold with a Poisson tensor compatible with
the Lorentzian metric in the sense of Hawkins and such that the
Poisson tensor cannot be constructed locally by commuting Killing
vector fields. This shows that Theorem 6.6 in [9] is false in the
Lorentzian case.

 For $\la=(\la_1,\la_2)\in\reel^2$, $0<\la_1\leq\la_2$, the
 oscillator group of dimension 6 is the connected and simply
 connected Lie group $G_\la$ whose Lie algebra is
 $$\G_\la:=vect\{e_{-1},e_0,e_1,e_e,\check{e}_1,\check{e}_2\}$$with
 brackets
 $$[e_{-1},e_j]=\la_j\check{e}_j,\qquad
 [e_j,\check{e}_j]=e_0,\qquad [e_{-1},\check{e}_j]=-\la_j e_j$$and
 the unspecified brackets are either zero or given by
 antisymmetry. These groups were introduced  by Medina (see [12]) as the only
 non commutative simply connected solvable Lie groups which have a
 bi-invariant Lorentzian metric.
  Their
 discrete co-compact subgroups where classified in [12].
 These groups have discrete co-compact subgroups if and only if
 the set $\{\la_1,\la_2\}$ generates a discrete subgroup of
 $(\reel,+)$. Let $\Lambda$ be a discrete co-compact subgroup of
 $G_\la$. The bi-invariant Lorentzian metric on $G_\la$ defines a Lorentzian
  metric on $P:=G_\la/\Lambda$ and we get an action of $G_\la$ on $P$
  by isometries. Consider now $r\in\wedge^2\G_\la$ given by
 $$r=e_0\wedge e_1+e_2\wedge \check{e}_1.$$It is easy to check
 that $r$ is a solution of the Yang-Baxter equation and $Imr$ is
 a nilpotent and hence an
  unimodular Lie algebra. According to Theorem 1.1, we get
  a Poisson tensor on $P$ which is compatible with the Lorentzian
  metric. But,  the Poisson tensor is
  non parallel with respect to the metric contravariant
  connection and hence it cannot be constructed locally from
  commuting Killing vectors fields.\end{enumerate}

{\bf Remark.} \begin{enumerate}\item In Theorem 1.2, if the action
of $Imr$ is not locally free, in general, the metacurvature does
not vanishe. For instance, consider the 2-dimensional Lie algebra
$\G=Vect\{e_1,e_2\}$ with $[e_1,e_2]=e_1$ and the action
$\Ga:\G\too\X(\reel^2)$ given by
$$\Ga(e_1)=\frac{\partial}{\partial x}\quad\mbox{and}\quad
\Ga(e_2)=x\frac{\partial}{\partial x}.$$ For $r=e_1\wedge e_2$,
$\Ga(r)$ is the trivial Poisson tensor and
$$\D^\Ga_\al\be=\al\left(\frac{\partial}{\partial x}\right)
\be\left(\frac{\partial}{\partial x}\right)dx.$$ A differential
1-form $\ga$ is parallel with respect $\D^\Ga$ if and only if
$\ga=f(x,y)dy$. In this case, for any $\al$ and $\be$, we have
$$\D^\Ga_\al\D^\Ga_\be d\ga=
\al\left(\frac{\partial}{\partial x}\right)
\be\left(\frac{\partial}{\partial x}\right) \frac{\partial
f}{\partial x}dx\wedge dy.$$This shows that the metacurvature does
not vanishe (see (7)).\end{enumerate}

 Section 2 is devoted to a complete proof of
Theorem 1.1 and Theorem 1.2.\bigskip

{\bf Notations.} For a smooth manifold $P$, $C^\infty(P)$ will
denote the space of smooth functions on $P$, $\Ga(P,V)$ will
denote the space of smooth sections  of a vector bundle,
$\Om^p(P):=\Ga(P,\wedge^pT^*P)$ and $\X^p(P):=\Ga(P,\wedge^pTP)$.
Lower case  Greek characters $\al,\be,\ga$ will mostly denote
1-forms. However, $\pi$ will denote a Poison bivector field and
$\om$ will denote a symplectic form.

For a manifold $P$ with a Poisson tensor $\pi$, $\an:T^*P\too TP$
will denote the anchor map given by $\be(\an(\al))=\pi(\al,\be),$
and $[\;,\;]_\pi$ will denote the Koszul  bracket given by
$$[\al,\be]_\pi=
 L_{\pi_{\#}(\al)}\be-L_{\pi_{\#}(\be)}\al-d(\pi(\al,\be)).$$
We will denote $P^{reg}$ the dense open set where the rank of
$\pi$ is locally constant.
 \section{Proof of Theorem 1.1 and Theorem 1.2}
 \subsection{Preliminaries}
Contravariant connections associated to a Poisson structure have
recently turned out to be useful in several areas of Poisson
geometry. Contravariant connections were defined by Vaismann [16]
and were analyzed in detail by Fernandes [7]. This notion appears
extensively  in the context of noncommutative deformations see
[8], [9], [15]. One can consult [7] for the general properties
 of contravariant connections.

In this subsection, I will give some general properties of
contravariant connections. Namely, I will recall the definition of
the metacurvature of a flat and torsion-free contravariant
connection $\D$, and I will give a necessary and sufficient
condition for the vanishing  of the metacurvature  in the case
where $\D$ is an $\F$-connection (see [7]).

Let $(P,\pi)$ be a Poisson manifold and $V\stackrel{p}\too P$ a
vector bundle over $P$. A contravariant connection on $V$ with
respect to $\pi$ is a map $\D:\Om^1(P)\times
\Ga(P,V)\too\Ga(P,V)$, $(\al,s)\mapsto\D_\al s$ satisfying the
following properties:\begin{enumerate} \item $\D_\al s$ is linear
over $C^\infty(P)$ in $\al$:
$$\D_{f\al+h\be}s=f\D_\al s+h\D_\be s,\quad f,g\in C^\infty(P);$$
\item $\D_\al s$ is linear over $\reel$ in $s$:
$$\D_{\al}(as_1+bs_2)=a\D_\al s_1+b\D_\al s_2,\quad a,b\in \reel;$$
\item $\D$ satisfies the following product rule:
$$\D_\al(fs)=f\D_\al s+\an(\al)(f)s,\quad f\in
C^\infty(P).$$\end{enumerate}

The curvature of a contravariant connection $\D$ is formally
identical to the usual definition
$$K(\al,\be)=\D_\al\D_\be-\D_\be\D_\al-\D_{[\al,\be]_\pi}.$$We
will call $\D$ flat if $K$ vanishes identically.

A contravariant connection $\D$ will be called an $\F$-connection
if its satisfies the following properties
$$\an(\al)=0\qquad\Rightarrow \qquad\D_\al=0.$$We will call $\D$
un $\F^{reg}$-connection if the restriction of $\D$ to $P^{reg}$
is an $\F$-connection.

If $V=T^*P$, one can define the torsion $T$ of $\D$ by
$$T(\al,\be)=\D_\al\be-\D_\be\al.$$

Let us define now an interesting  class of contravariant
connections, namely contravariant connection associated naturally
to a Poisson tensor and a pseudo-Riemannian metric.

 Let $P$ be a pseudo-Rimannian
manifold and $\pi$ a Poisson tensor on $P$.

The metric contravariant connection associated naturally to
$(\pi,<\;,\;>)$ is the unique contravariant connection $D$ such
that: \begin{enumerate}\item the  metric $<,>$ is parallel with
respect to $D$, i.e.,
$$\pi_{\#}(\al).<\be,\ga>=<D_\al\be,\ga>+<\be,D_\al\ga>;$$
\item $D$ is torsion-free.
\end{enumerate}
 One can define $D$ by the Koszul formula
\begin{eqnarray*}
2<D_\al\be,\ga>&=&\pi_{\#}(\al).<\be,\ga>+\pi_{\#}(\be).<\al,\ga>-
\pi_{\#}(\ga).<\al,\be>\\
&+&<[\ga,\al]_\pi,\be>+<[\ga,\be]_\pi,\al>+<[\al,\be]_\pi,\ga>.
\qquad(4)\end{eqnarray*}

Let us recall briefly the definition of the metacurvature. For
details see [9].

Let $(P,\pi)$ be a Poisson manifold and $\D$ a torsion-free and
flat contravariant connection with respect to $\pi$. Then, there
exists a bracket $\{\;,\;\}$ on the differential graded algebra of
differential forms $\Om^*(P)$ such that: \begin{enumerate} \item
$\{\;,\;\}$ is $\reel$-bilinear, degree 0 and antisymmetric, i.e.
$$\{\sigma,\rho\}=-(-1)^{deg\sigma deg\rho}\{\rho,\sigma\}.
$$
\item The differential $d$ is a derivation with respect to
$\{\;,\;\}$ i.e.
$$d\{\sigma,\rho\}=\{d\sigma,\rho\}+(-1)^{deg\sigma}\{\sigma,d\rho\}.
$$ \item $\{\;,\;\}$ satisfies the product rule
$$\{\sigma,\rho\wedge\la\}=\{\sigma,\rho\}\wedge\la+(-1)^{deg\sigma
deg\rho}\rho\wedge\{\sigma,\la\}.$$ \item For any $f,h\in
C^\infty(P)$ and for any $\sigma\in\Om^*(P)$ the bracket $\{f,g\}$
coincides with the initial Poisson bracket and
$$\{f,\sigma\}=\D_{df}\sigma.$$\end{enumerate}
Hawkins called this bracket a generalized Poisson bracket and
showed that there exists a $(2,3)$-tensor $\M$ such that the
following assertions are equivalent:\begin{enumerate}\item The
generalized Poisson bracket satisfies the graded Jacobi identity
$$\{\{\sigma,\rho\},\la\}=\{\sigma,\{\rho,\la\}\}
-(-1)^{deg\sigma deg\rho}\{\rho,\{\sigma,\la\}\}.
$$\item The tensor $\M$ vanishes identically.\end{enumerate}
$\M$ is called the metacurvature and is given by
$$\M(df,\al,\be)=\{f,\{\al,\be\}\}-\{\{f,\al\},\be\}-
\{\{f,\be\},\al\}.\eqno(5)$$ Hawkins pointed out in $[9]$ pp. 9,
that for any parallel 1-form $\al$ and any 1-form $\be$, the
generalized Poisson bracket of $\al$ and $\be$ is given by
$$\{\al,\be\}=-\D_{\be}d\al.\eqno(6)$$Then, one can deduce from
$(5)$ that for any parallel 1-form $\al$ and for any $\be,\ga$, we
have
$$\M(\al,\be,\ga)=-\D_\be\D_\ga d\al.\eqno(7)$$

The definition of a contravariant connection is similar to the
definition of an ordinary (covariant) connection, except that
cotangent vectors have taken the place of tangent vectors. So one
can translate many definitions, identities and proof for covariant
connections to contravariant connections simply by exchanging the
roles of tangent and cotangent vectors and replacing Lie Bracket
with Koszul bracket. Nevertheless, there are some differences
between those two notions.  Fernandes pointed out in [7] that the
equation $D\al=0$ cannot be solved locally for a general flat
contravariant connection $D$. However, he showed that for a flat
$\F$-connection this equation  can be solved locally. We will give
now a proof of this fact which is different from Fernandes's
proof.
\begin{pr} Let $(P,\pi)$ be a Poisson manifold and $\D$ a
flat contravariant connection with respect $\pi$. Let $p$ be a
regular point of $\pi$ and
$(p_1,\ldots,p_r,q_1,\ldots,q_r,z_1,\ldots,z_l)$ a Darboux
coordinates on a neighborhood $\U$ of $p$. We denote by $N$ the
submanifold given by $z_i=0$, $i=1,\ldots,l$. Suppose that the
restriction of $\D$ to $\U$ is an $\F$-connection. Then, for any
$\be\in \Ga(N,T_{|N}^*P)$ there exists a unique 1-form
$\wi\be\in\Om^1(\U)$ such that $\D\wi\be=0$ and
$\wi\be_{|N}=\be$.\end{pr}

{\bf Proof.} In the Darboux coordinates system
$(p_1,\ldots,p_r,q_1,\ldots,q_r,z_1,\ldots,z_l)$, we have
$$\pi=\sum_{i=1}^r\frac{\partial }{\partial p_i}\wedge\frac{\partial }{\partial
q_i}.$$ The connection $\D$ is entirely determined by the
Christoffel symbols
$\Ga_{p_ip_j}^{p_k},\Ga_{p_ip_j}^{q_k},\Ga_{p_ip_j}^{z_k}$ and so
on. We are looking for
$$\wi\be=\sum_{i=1}^r(a_idp_i+b_idq_i)+\sum_{i=1}^lc_idz_i$$ such
that  $\D\wi\be=0$ and $\wi\be$ along $N$ coincides with $\be$.

Since $\D$ is a torsion-free $\F$-connection, and since $dz_i$,
$i=1,\ldots,l$, is in the center of $(\Om^1(\U),[\;,\;]_\pi)$, one
have $\D_{dz_i}=Ddz_i=0$. Hence $\D\wi\be=0$ is equivalent to the
following two systems:
$$\left\{\begin{array}{lll}
\di\frac{\partial a_k}{\partial
q_i}&=&-\di\sum_{l=1}^r(a_l\Ga_{p_ip_l}^{p_k}+b_l\Ga_{p_iq_l}^{p_k}),\\
\di\frac{\partial a_k}{\partial
p_i}&=&\di\sum_{l=1}^r(a_l\Ga_{q_ip_l}^{p_k}+b_l\Ga_{q_iq_l}^{p_k}),\\
\di\frac{\partial b_k}{\partial
q_i}&=&-\di\sum_{l=1}^r(a_l\Ga_{p_ip_l}^{q_k}+b_l\Ga_{p_iq_l}^{q_k}),\\
\di\frac{\partial b_k}{\partial
q_i}&=&\di\sum_{l=1}^r(a_l\Ga_{q_ip_l}^{q_k}+b_l\Ga_{q_iq_l}^{q_k}).
\end{array}\right.\eqno(*)$$
$$\left\{\begin{array}{lll}\di
\frac{\partial c_k}{\partial
q_i}&=&-\di\sum_{l=1}^r(a_l\Ga_{p_ip_l}^{z_k}+b_l\Ga_{p_iq_l}^{z_k}),\\
\di\frac{\partial c_k}{\partial
p_i}&=&\di\sum_{l=1}^r(a_l\Ga_{q_ip_l}^{z_k}+b_l\Ga_{q_iq_l}^{z_k}).
\end{array}\right.\eqno(**)$$
One can see the functions $(a_1,\ldots,a_r,b_1,\ldots,b_r)$ in
$(*)$ as functions with  variables $(q_i,p_i)$ and parameters
$(z_1,\ldots,z_l)$. The vanishing of curvature gives necessary
integrability conditions of $(*)$ and hence for any initial value
and any value of the parameter there exists a unique solution
(which depends smoothly on the parameter)
$(a_1,\ldots,a_r,b_1,\ldots,b_r)$ of $(*)$. For $k=1,\ldots,l$,
consider the 1-form with  variables $(q_i,p_i)$ and parameters
$(z_1,\ldots,z_l)$
$$\al_k=\sum_{i=1}^r
\left(-\left(\sum_{l=1}^r(a_l\Ga_{p_ip_l}^{z_k}+b_l\Ga_{p_iq_l}^{z_k})
\right)dq_i+\left(\sum_{l=1}^r(a_l\Ga_{q_ip_l}^{z_k}+b_l\Ga_{q_iq_l}^{z_k})\right)dp_i
\right).$$The vanishing of the curvature implies that $d\al_k=0$
and hence there exists a function $c_k$ such that $dc_k=\al_k$.
This solves $(**)$. $\Box$

By combining $(7)$ and Proposition 2.1, we get the following
useful Proposition. \begin{pr} Let $(P,\pi)$ be a Poisson manifold
and $\D$ a torsion-free and flat $\F^{reg}$-connection with
respect to $\pi$. Then the metacurvature of $\D$ vanishes if and
only if, for any local parallel 1-form on $P^{reg}$,
$D^2d\al=0$.\end{pr}
\subsection{Proof of Theorem 1.1}
Let $P$ be a differentiable manifold, $\G$ a Lie algebra with
$\Ga:\G\too \X(P)$ a Lie algebras morphism from $\G$ to the Lie
algebra of  vector fields and $r\in\wedge^2\G$ a solution of the
classical Yang-Baxter equation such that $Imr$ is an unimodular
Lie algebra. There exists a basis
$(e_1,\ldots,e_n,f_1,\ldots,f_n)$ of $Imr$ such that symplectic
form $\om_r$ is given by
$$\om_r=\sum_{i=1}^ne_i^*\wedge f_i^*.$$ Since $Imr$ is
unimodular, then for any $z\in Imr$, the trace of $ad_z$ is null.
This is equivalent to
$$\sum_{i=1}^n\om_r([z,e_i],f_i)+\om_r(e_i,[z,f_i])=0.$$
From $(1)$, one get that this relation is equivalent to
$$\sum_{i=1}^n\om_r(z,[e_i,f_i])=0$$and hence to
$$\sum_{i=1}^n[e_i,f_i]=0.$$Now let $\e$ a volume form on $P$ such
that $L_{\Ga(e_i)}\e=L_{\Ga(f_i)}\e=0$ for $i=1,\ldots,n$. We have
\begin{eqnarray*}
d(i_{\Ga(r)}\e)&=&d\left(\sum_{i=1}^ni_{\Ga(e_i)\wedge\Ga(f_i)}\e\right)\\
&=&\sum_{i=1}^n\left(i_{[\Ga(e_i),\Ga(f_i)]}\e-
i_{\Ga(e_i)}L_{\Ga(f_i)}\e-i_{\Ga(f_i)}L_{\Ga(e_i)}\e\right)\\
&=&i_{\Ga(\sum_{i=1}^n[e_i,f_i])}\e=0.\end{eqnarray*} This gives a
proof of Theorem 1.1.$\Box$ \subsection{Proof of Theorem 1.2} Let
$P$ be a differentiable manifold, $\G$ a Lie algebra with
$\Ga:\G\too \X(P)$ a Lie algebras morphism from $\G$ to the Lie
algebra of vector fields and $r\in\wedge^2\G$ a solution of the
classical Yang-Baxter equation. If a basis $(u_1,\ldots,u_n)$  of
$\G$ is chosen then we can write $r=\sum_{i<j}a_{ij}u_i\wedge
u_j$. For any $\al,\be$ two 1-forms on $P$, we put
$$\D^{\Ga}_\al\be=\sum_{i<j}a_{ij}\left(\al(U_i)L_{U_j}\be-\al(U_j)L_{U_i}\be
\right),$$where $U_i=\Ga(u_i)$. One can check easily that this
formula defines a contravariant connection with respect to
$\Ga(r)$ which depends only on $r$ and $\Ga$ and doesn't depend on
the basis $(u_1,\ldots,u_n)$. One can also check that $\D^\Ga$ is
torsion-free. If, for any $u\in Imr$, $\Ga(u)$ is a Killing vector
field, then $\D^\Ga$ is the metric contravariant connection
associated to the metric and the Poisson tensor $\Ga(r)$.

Let us now compute the  curvature of $\D^\Ga$ and show that it
vanishes identically.

There exists    a basis $(u_1,\ldots,u_p,v_1,\ldots,v_p)$ of $Imr$
such that $r=\sum_{i=1}^pu_i\wedge v_i.$ We denote $U_i=\Ga(u_i)$
and $V_i=\Ga(v_i)$ for $i=1,\ldots,p$. We have
$$\D^\Ga_\al\be=L_{\an(\al)}\be+\sum_{i=1}^pA^i(\al,\be)$$where
$A^i(\al,\be)=\be(U_i)d\left(\al(V_i)\right)-\be(V_i)d\left(\al(U_i)\right)$
and $\an$ is the anchor map associated to $\Ga(r)$. With this in
mind, we get for any $f,g,h\in C^\infty(P)$,
\begin{eqnarray*}
K(df,dg,dh)&=&\sum_{i=1}^p\left(A^i(df,d\{g,h\})-A^i(dg,d\{f,h\})\right.\\
&+&\left.
L_{\an(df)}A^i(dg,dh)-L_{\an(dg)}A^i(df,dh)-A^i(d\{f,g\},dh)\right)\\
&+&\sum_{i,j=1}^p\left(A^j(df,A^i(dg,dh))-A^j(dg,A^i(df,dh)).\right)
\end{eqnarray*}

A straightforward computation gives{\small
\begin{eqnarray*}
K(df,dg,dh)&&=\\\sum_{i,j=1}^p&\left\{\right.&
\di\left(U_j(g)[U_i,V_j](h)-U_j(h)[U_i,V_j](g)+
V_j(g)[U_j,U_i](h)-V_j(h)[U_j,U_i](g)\right)d(V_i(f))\\
&+&\left(U_j(g)[V_j,V_i](h)-U_j(h)[V_j,V_i](g)+V_j(g)[V_i,U_j](h)-
V_j(h)[V_i,U_j](g)\right)d(U_i(f))\\
&-&\left(U_j(f)[U_i,V_j](h)-U_j(h)[U_i,V_j](f)+
V_j(f)[U_j,U_i](h)-V_j(h)[U_j,U_i](f)\right)d(V_i(g))\\
&-&\left(U_j(f)[V_j,V_i](h)-U_j(h)[V_j,V_i](f)+V_j(f)[V_i,U_j](h)-V_j(h)
[V_i,U_j](f)\right)d(U_i(g))\\
&+&U_i(h)V_j(g)d\left([U_j,V_i](f)\right)-U_i(h)V_j(f)d\left([U_j,V_i](g)\right)\\
&+&U_i(h)U_j(g)d\left([V_i,V_j](f)\right)-U_i(h)U_j(f)d\left([V_i,V_j](g)\right)\\
&+&V_i(h)U_j(g)d\left([V_j,U_i](f)\right)-V_i(h)U_j(f)d\left([V_j,U_i](g)\right)\\
&+&\left.V_i(h)V_j(g)d\left([U_i,U_j](f)\right)-V_i(h)V_j(f)d\left([U_i,U_j](g)
\right)\right\}.\end{eqnarray*}}

The vanishing of $K$ is a consequence of the equation $[r,r]=0$
which is equivalent to
$$\om_r([x,y],z)+\om_r([y,z],x)+\om_r([z,x],y)=0\quad\forall x,y,z\in Imr.
\eqno(*)$$

Now $\om_r=\sum_{i=1}^pu_i^*\wedge v_i^*$ where
$(u_1^*,\ldots,u_p^*,v_1^*,\ldots,v_p^*)$ is the dual basis of
$(u_1,\ldots,u_p,v_1,\ldots,v_p)$. If one write
$[u_i,u_j]=\sum_{k=1}^p(C_{u_iu_j}^{u_k}u_k+C_{u_iu_j}^{v_k}v_k)$
and so on, one can see easily that the condition $(*)$ is
equivalent to
$$\left\{\begin{array}{ccc}
C_{v_jv_k}^{u_i}+C_{v_kv_i}^{u_j}+C_{v_iv_j}^{u_k}&=&0,\\
C_{u_ju_k}^{v_i}+C_{u_ku_i}^{v_j}+C_{u_iu_j}^{v_k}&=&0,\\
C_{u_ju_k}^{u_i}-C_{u_kv_i}^{v_j}-C_{v_iu_j}^{v_k}&=&0,\\
C_{v_iv_j}^{v_k}-C_{u_kv_i}^{u_j}-C_{v_ju_k}^{u_i}&=&0.\end{array}\right.\qquad\forall
i,j,k.$$ For $i=1,\ldots,p$, the coefficients of $d(V_i(f))$ and
$d(U_i(f))$ in the expression of $K(df,dg,dh)$ are respectively
\begin{eqnarray*}
\sum_{j,k=1}^p\left(C_{u_iv_j}^{u_k}-C_{u_iv_k}^{u_j}+C_{v_kv_j}^{v_i}\right)U_j(g)U_k(h)
+\left(C_{u_iv_j}^{v_k}-C_{u_ku_i}^{u_j}+C_{v_ju_k}^{v_i}\right)U_j(g)V_k(h)\\
+\left(-C_{u_iv_j}^{v_k}+C_{u_ku_i}^{u_j}+C_{u_kv_j}^{v_i}\right)U_j(h)V_k(g)+
\left(C_{u_ju_i}^{v_k}-C_{u_ku_i}^{v_j}+C_{u_ku_j}^{v_i}\right)V_j(g)V_k(h),\end{eqnarray*}
\begin{eqnarray*}
\sum_{j,k=1}^p\left(C_{v_jv_i}^{u_k}-C_{v_kv_i}^{u_j}+C_{v_kv_j}^{u_i}\right)U_j(g)U_k(h)
+\left(C_{v_jv_i}^{v_k}-C_{v_iu_k}^{u_j}+C_{v_ju_k}^{u_i}\right)U_j(g)V_k(h)\\
+\left(-C_{v_jv_i}^{v_k}+C_{v_iu_k}^{u_j}+C_{u_kv_j}^{u_i}\right)U_j(h)V_k(g)+
\left(C_{v_iu_j}^{v_k}-C_{v_iu_k}^{v_j}+C_{u_ku_j}^{u_i}\right)V_j(g)V_k(h).\end{eqnarray*}
Those coefficients vanish according to the relations above. The
same thing will happens for the coefficients of $d(V_i(g))$ and
$d(U_i(g)$. This shows that $K$ vanishes identically.

Suppose now that the action of $Imr$ on $P$ is locally free. This
will implies obviously that $\D^\Ga$ is an $\F$-connection and a
1-form $\be$ is parallel with respect $\D^\Ga$ if and only if
$L_{\Ga(u)}\be=0$ for all $u\in Imr$. For any parallel 1-form
$\be$, we have $L_{\Ga(u)}d\be=0$ and hence $\D^\Ga d\be=0$.
According to Proposition 2.2, this implies that the metacurvature
of $\D^\Ga$ vanishes identically, which achieves the proof of
Theorem 1.2.$\Box$

\eject

{\bf References}\bigskip

[1] {\bf Benson G.,  Gordon C.,} {\it K\"ahler and symplectic
structures on nilmanifolds,} Topology {\bf 27,4} (1988), 513-518.

[2] {\bf Bordemann M., Medina A., Ouadfel A.,} {\it Le groupe
affine comme vari\'et\'e symplectique}, T\^{o}hoku Math. J., {\bf
45} (1993), 423-436.

 [3] {\bf
Boucetta M., } {\it Compatibilit\'e des structures
pseudo-riemanniennes et des structures de Poisson}{ C. R. Acad.
Sci. Paris, {\bf t. 333}, S\'erie I, (2001) 763--768.}

[4] {\bf Boucetta M., } {\it Poisson manifolds with compatible
pseudo-metric and pseudo-Riemannian Lie algebras}{ Differential
Geometry and its Applications, {\bf Vol. 20, Issue 3}(2004),
279--291.}

[5] {\bf Chu B. Y.,} {\it Symplectic homogeneous spaces,} Trans.
Am. Math. Soc. {\bf 197} (1974), 145-159.

[6] {\bf  De Smedt V.,} {\it Existence of a Lie bialgebra
structure on every Lie algebra,} Letters in Mathematical Physics,
{\bf 31} (1994), 225-231.

 [7] {\bf Fernandes R. L.,} {\it Connections in
Poisson Geometry1: Holonomy and invariants}, J. of Diff. Geometry
{\bf 54} (2000), 303-366.

[8] {\bf Hawkins E.,} {\it Noncommutative Rigidity}, Commun. Math.
Phys. {\bf 246} (2004) 211-235. math.QA/0211203.

[9] {\bf Hawkins E.,} {\it The structure of noncommutative
deformations,}

arXiv:math.QA/0504232.

 [10] {\bf A. Lichnerowicz and A. Medina,} {\it On Lie groups with
left-invariant symplectic or K\"ahlerian structures,} {Letters in
Mathematical Physics {\bf16} (1988), 225--235.}

[11] {\bf  Medina A.,  Revoy Ph.,} {\it Alg\`ebres de Lie et
produit scalaire invariant,} {Ann. Ec. Norm. Sup., 4\`eme s\'erie,
{\bf t. 18} (1985), 553--561}.

[12] {\bf  Medina A.,  Revoy Ph.,} {\it Les groupes oscillateurs
et leurs r\'eseaux,} Manuscripta Math. {\bf 52} (1985), 81-95.

[13] {\bf  Medina A.,  Revoy Ph.,} {\it Groupes de Lie \`a
structure symplectique invariante,} Symplectic geometry, groupoids
and integrable systems, in "S\'eminaire Sud Rhodanien", M.S.R.I.,
New York/Berlin: Springer-Verlag, (1991), 247-266.

[14] {\bf Milnor J.,} {\it Curvature of left invariant metrics on
Lie groups,} Adv. in Math. {\bf21} (1976), 283-329.

[15] {\bf Reshetikhin N., Voronov A., Weinstein A.,} {\it
Semiquantum Geometry, Algebraic geometry}, J. Math. Sci. {\bf
82(1)} (1996) 3255-3267.

 [16] {\bf Vaisman I.}, {\it Lecture on
the geometry of Poisson manifolds}, Progr. In Math. {\bf Vol.
118}, Birkhausser, Berlin, (1994).

\vskip3cm

M. Boucetta\\
Facult\'e des Sciences et Techniques \\
BP 549 Marrakech\\
Morocco
\\
Email: {\it boucetta@fstg-marrakech.ac.ma

mboucetta2@yahoo.fr}

\end{document}